\documentclass{article}
\usepackage{graphicx,bm,url,amssymb,amsmath,hyperref,color,geometry}

\def\no{\noindent}
\def\pmatrix{\left(\begin{array}}
\def\endpmatrix{\end{array}\right)}

\def\RR{\mathbb{R}}

\def\dd{\mathrm{d}}

\def\aa{\alpha}

\def\cpr{$^\copyright\,$}

\begin{document}

\title{Solving FDE-IVPs by using Fractional HBVMs: \\ some experiments with the {\tt fhbvm} code}

\author{Luigi Brugnano\,\footnote{
 Dipartimento di Matematica e Informatica ``U.\,Dini'', 
             Universit\`a di Firenze,  Italy\\    \url{luigi.brugnano@unifi.it} \quad
                       \url{https://orcid.org/0000-0002-6290-4107}}
   \and  Gianmarco Gurioli\,\footnote{
              Dipartimento di Matematica e Informatica ``U.\,Dini'',
             Universit\`a di Firenze,  Italy\\
              \url{gianmarco.gurioli@unifi.it}  \quad
              \url{https://orcid.org/0000-0003-0922-8119}}
     \and Felice Iavernaro\,\footnote{         Dipartimento di Matematica, 
           Universit\`a di Bari Aldo Moro, Italy\\
            \url{felice.iavernaro@uniba.it}\quad
            \url{https://orcid.org/0000-0002-9716-7370}}}

\maketitle

\begin{abstract} In this paper we report a few numerical tests by using a slight extension of the Matlab\cpr code {\tt fhbvm} in \cite{BGI2024}, implementing {\em Fractional HBVMs}, a recently introduced class of numerical methods for solving Initial Value Problems of Fractional Differential Equations (FDE-IVPs).  The reported experiments are aimed to give evidence of its effectiveness.

\medskip
\no{\bf Keywords:} fractional differential equations, fractional integrals, Caputo derivative, Jacobi polynomials, Fractional HBVMs, FHBVMs.

\medskip
\no{\bf MSC:} 34A08, 65R20, 65-04.
\end{abstract}

\section{Introduction} 

Fractional differential equations have become popular in many applications (see, e.g.,  the classical references \cite{Di2010,Po1999} for an introduction). For this reason, their numerical solution is an active field of investigation.
In particular, we here report a few numerical tests obtained by using the Matlab\cpr code {\tt fhbvm} \cite{BGI2024},
available at the URL \cite{fhbvm}. The code itself is based on a class of methods named {\em Fractional HBVMs (FHBVMs)}, recently introduced in \cite{BBBI2023}, as a generalization of {\em Hamiltonian Boundary Value Methods (HBVMs)}, a class of low-rank Runge-Kutta methods initally devised for the efficient numerical  solution of Hamiltonian dynamical systems (see, e.g., \cite{BI2016,BI2018}), and, subsequently, generalized along several directions (see, e.g., \cite{ABI2019,ABI2022-1,ABI2023,BBBI2023,BFCIV2022,BI2022}). A main feature of HBVMs is the fact that they can gain spectral accuracy when approximating ODE-IVPs \cite{ABI2020,BMIR2019,BMR2019}, and this holds true also in the FDE case \cite{BBBI2023}. In this paper we shall use a slight extension of the code in \cite{BGI2024}, modified for handling  FDE-IVPs in the form:\footnote{The code in \cite{BGI2024} was conceived for the case $\ell=1$.}
\begin{eqnarray}\nonumber
y^{(\aa)}(t) &=& f(y(t)), \qquad t\in[0,T],\qquad y^{(i)}(0) = y_0^i\in\RR^m, \quad i=0,\dots,\ell-1, \\
\aa &\in&(\ell-1,\ell), \label{ivp}
\end{eqnarray}
where ~$y^{(\aa)}(t)$~ is the Caputo fractional derivative and, for the sake of brevity, we have omitted the argument $t$ for $f$. It is known that under suitable regularity assumptions (see, e.g. \cite{DGS2020,Ga2018}), the solution of (\ref{ivp}) is given by: 
\begin{equation}\label{sol0}
y(t) = \sum_{i=0}^{\ell-1}\frac{t^i}{i!}y_0^i + I^{\aa} f(y(t)) \equiv  \sum_{i=0}^{\ell-1}\frac{t^i}{i!}y_0^i +  \frac{1}{\Gamma(\aa)}\int_0^t (t-x)^{\aa-1} f(y(x))\dd x, \qquad t\in[0,T],
\end{equation}
where $I^\aa f(y(t))$ is the Riemann-Liouville integral of $f(y(t))$. With this premise, the structure of the paper is as follows: in Section~\ref{start} we recall the main facts about the numerical solution of FDE-IVPs proposed in \cite{BBBI2023};  in Section~\ref{num} we report some numerical experiments with the code, also providing some comparisons with another existing one, to give evidence of its effectiveness; at last, a few conclusions are given in Section~\ref{fine}.

\section{Fractional HBVMs}\label{start}
To begin with, in order to obtain a piecewise approximation to the solution of the problem,  we consider a partition of the integration interval in the form:
\begin{equation}\label{tn}
t_n = t_{n-1} + h_n,  \qquad h_n>0, \qquad n=1,\dots,N, \qquad t_0=0, \quad t_N=T. 
\end{equation}
In particular, the code {\tt fhbvm} implements either a {\em graded mesh},
\begin{equation}\label{hnr}
h_n = r h_{n-1} \equiv r^{n-1} h_1, \qquad n=1\dots,N,
\end{equation}
where $r>1$ and $h_1>0$ satisfy, by virtue of (\ref{tn})--(\ref{hnr}),
~$
h_1\dfrac{r^N-1}{r-1} = T,
$~
or a {\em uniform mesh},
\begin{equation}\label{unif}
h_n \equiv h_1 := \frac{T}N, \quad n=1,\dots,N, \quad \Rightarrow\quad t_n = nh_1, \quad n=0,\dots,N.
\end{equation}
The former mesh is clearly appropriate when the vector field in (\ref{ivp}) is not very smooth at the origin, whereas the latter mesh is recommended, otherwise. The code automatically selects the most appropriate mesh and its parameters \cite{BGI2024}.
The method is aimed at obtaining a piecewise approximation $\sigma(t)$ defined over the partition (\ref{tn}). 
In order to derive it, let us set,  
\begin{equation}\label{yn}
y_n(ch_n) := y(t_{n-1}+ch_n), \quad \sigma_n(ch_n) := \sigma(t_{n-1}+ch_n), \quad c\in[0,1],\qquad n=1,\dots,N,
\end{equation}
the restriction of the solution of (\ref{ivp}) and of its approximation $\sigma(t)$ on the interval $[t_{n-1},t_n]$, respectively. Consequently, the problem (\ref{ivp}) can be rewritten as:
\begin{equation}\label{ivpN}
\left\{\begin{array}{rcl}
y_n^{(\aa)}(ch_n) &=& f(y_n(ch_n)), \qquad c\in[0,1], \qquad n=1,\dots,N,\\[2mm]
y_1^{(i)}(0) &=& y_0^i, \qquad i=0,\dots,\ell-1. \end{array}\right.
\end{equation}
It is to be noticed that the r.h.s. of each local problem in (\ref{ivpN}) only depends on the corresponding restriction of the solution, i.e., $y_n$. Further, we consider the expansion of the local vector fields along the following basis of orthonormal Jacobi polynomials:\footnote{Here, $\aa>0$ is the same parameter as in (\ref{ivp}).} 
\begin{equation}\label{jaco}
P_i\in\Pi_i, \qquad \int_0^1 \omega(c)P_i(c)P_j(c)\dd c = \delta_{ij}, \qquad i,j=0,1,\dots,\qquad \omega(c) := \aa(1-c)^{\aa-1}.\end{equation}
Consequently, one obtains:
\begin{equation}\label{exp1}
f(y_n(ch_n)) = \sum_{j\ge0} P_j(c)\gamma_j(y_n), \qquad c\in[0,1],\qquad \gamma_j(y_n) = \int_0^1\omega(\tau)P_j(\tau)f(y_n(\tau h_n))\dd\tau,
\end{equation}
so that (\ref{ivpN}) can be equivalently rewritten as
\begin{equation}\label{ivpN1}
\left\{\begin{array}{rcl}
y_n^{(\aa)}(ch_n) &=& \sum_{j\ge0} P_j(c)\gamma_j(y_n), \qquad c\in[0,1], \qquad n=1,\dots,N,\\[2mm]
y_1^{(i)}(0) &=& y_0^i, \qquad i=0,\dots,\ell-1. \end{array}\right.
\end{equation}
The approximation $\sigma(t)$ is then defined by truncating the above series to finite sums with $s\ge1$ terms:
\begin{equation}\label{ivpN2}
\left\{\begin{array}{rcl}
\sigma_n^{(\aa)}(ch_n) &=& \sum_{j=0}^{s-1} P_j(c)\gamma_j(\sigma_n), \qquad c\in[0,1], \qquad n=1,\dots,N,\\[2mm]
\sigma_1^{(i)}(0) &=& y_0^i, \qquad i=0,\dots,\ell-1, \end{array}\right.
\end{equation}
with $\gamma_j(\sigma_n)$ defined according to (\ref{exp1}), by formally replacing $y_n$ with $\sigma_n$.

\medskip
For the continuous solution, for $t\,\equiv\, t_{n-1}+ch_n$, $c\in[0,1]$, one then obtains (see (\ref{tn})):
\begin{eqnarray*}
y(t) &\equiv&y_n(ch_n) ~=~ \sum_{i=0}^{\ell-1}\frac{t^i}{i!}y_0^i ~+~ I^\aa f(y(t_{n-1}+ch_n))\\[1mm]
&=& \sum_{i=0}^{\ell-1}\frac{t^i}{i!}y_0^i  ~+~ \frac{1}{\Gamma(\aa)}\int_0^{t_{n-1}+ch_n} (t_{n-1}+ch_n-x)^{\aa-1} f(y(x))\dd x\\[2mm]
&=&\sum_{i=0}^{\ell-1}\frac{t^i}{i!}y_0^i  ~+~ \frac{1}{\Gamma(\aa)}\sum_{\nu=1}^{n-1} \int_{t_{\nu-1}}^{t_\nu} (t_{n-1}+ch_n-x)^{\aa-1}f(y(x))\dd x\\[2mm]  
&&\qquad + ~\frac{1}{\Gamma(\aa)} \int_{t_{n-1}}^{t_{n-1}+ch_n} (t_{n-1}+ch_n-x)^{\aa-1}f(y(x))\dd x\\[2mm]  
&=&\sum_{i=0}^{\ell-1}\frac{t^i}{i!}y_0^i  ~+~ \frac{1}{\Gamma(\aa)}\sum_{\nu=1}^{n-1} \int_0^{h_\nu} (t_{n-1}-t_{\nu-1}+ch_n-x)^{\aa-1}f(y_\nu(x))\dd x\\[2mm]  
&&\qquad + ~\frac{1}{\Gamma(\aa)} \int_0^{ch_n} (ch_n-x)^{\aa-1}f(y_n(x))\dd x\\[2mm]  
&=&\sum_{i=0}^{\ell-1}\frac{t^i}{i!}y_0^i  ~+~ \frac{1}{\Gamma(\aa)}\sum_{\nu=1}^{n-1} 
h_\nu^\aa \int_0^1 \left( \frac{t_{n-1}-t_{\nu-1}}{h_\nu}+c\frac{h_n}{h_\nu}-\tau \right)^{\aa-1}f(y_\nu(\tau h_\nu))\dd \tau\\[2mm]  
&&\qquad + ~\frac{h_n^\aa}{\Gamma(\aa)} \int_0^c (c-\tau)^{\aa-1}f(y_n(\tau h_n))\dd \tau, \qquad c\in[0,1].
\end{eqnarray*}

\medskip
\no Consequently, by considering, hereafter, the case of the graded mesh (\ref{hnr}), and taking into account (\ref{ivpN1}), one has:\footnote{Similar arguments can be used in the case of the uniform mesh (\ref{unif}), see \cite{BBBI2023,BGI2024}.} 
\begin{eqnarray}\nonumber
y_n(ch_n)
&=&\sum_{i=0}^{\ell-1}\frac{t^i}{i!}y_0^i ~+~ \frac{1}{\Gamma(\aa)}\sum_{\nu=1}^{n-1} 
h_\nu^\aa \int_0^1 \left( \frac{r^{n-\nu}-1}{r-1}+cr^{n-\nu}-\tau \right)^{\aa-1}\sum_{j\ge0}P_j(\tau)\gamma_j(y_\nu)\dd\tau
\\[2mm]  \label{ynch}
&&\qquad + ~\frac{h_n^\aa}{\Gamma(\aa)} \int_0^c (c-\tau)^{\aa-1}\sum_{j\ge0}P_j(\tau)\gamma_j(y_n)\dd\tau,
\qquad c\in[0,1].
\end{eqnarray}
According to (\ref{ivpN2}), the corresponding restriction of the approximation $\sigma(t)$ on the interval $[t_{n-1},t_n]$  is obtained by truncating the infinite series in  (\ref{ynch}) to finite sums:
\begin{eqnarray}\nonumber
\sigma_n(ch_n)  
&=&\sum_{i=0}^{\ell-1}\frac{t^i}{i!}y_0^i ~+~ \frac{1}{\Gamma(\aa)}\sum_{\nu=1}^{n-1} 
h_\nu^\aa \int_0^1 \left( \frac{r^{n-\nu}-1}{r-1}+cr^{n-\nu}-\tau \right)^{\aa-1}\sum_{j=0}^{s-1}P_j(\tau)\gamma_j(\sigma_\nu)\dd\tau
\\[2mm]  \nonumber
&&\qquad + ~\frac{h_n^\aa}{\Gamma(\aa)} \int_0^c (c-\tau)^{\aa-1}\sum_{j=0}^{s-1}P_j(\tau)\gamma_j(\sigma_n)\dd\tau\\[1mm] \nonumber
&=:& \phi_{n-1}^{\aa,s}(c)  ~+ ~\frac{h_n^\aa}{\Gamma(\aa)} \int_0^c (c-\tau)^{\aa-1}\sum_{j=0}^{s-1}P_j(\tau)\gamma_j(\sigma_n)\dd\tau \\  \label{signch}
&\equiv& \phi_{n-1}^{\aa,s}(c)  ~+ h_n^\aa\sum_{j=0}^{s-1} I^\aa P_j(c)\gamma_j(\sigma_n),
\qquad c\in[0,1],
\end{eqnarray}
with $I^\aa P_j(c)$ the Riemann-Liouville integral of $P_j(c)$ (see  (\ref{sol0})). We observe that in (\ref{signch}) the {\em memory term} ~$\phi_{n-1}^{\aa,s}(c)$~ only depends on past values of the approximation (in fact, from $\sigma_\nu$, $\nu=1,\dots,n-1$). Consequently, the local problem consists in computing the Fourier coefficients $\gamma_j(\sigma_n)$, $j=0,\dots,s-1$. These latter coefficients can be approximated to full machine accuracy by using the Gauss-Jacobi formula of order $2k$ based at the zeros  $c_1,\dots,c_k$, of $P_k(c)$, with corresponding weights (see (\ref{jaco}))
$$b_i = \int_0^1 \omega(c)\ell_i(c)\dd c, \qquad \ell_i(c) = \prod_{j\ne i}\frac{c-c_j}{c_i-c_j}, \qquad i=1,\dots,k,$$ 
by choosing  a value of $k$, $k\ge s$, large enough. In other words,
\begin{equation}\label{gammajn}
\gamma_j^n \,:=\, \sum_{i=1}^k b_i P_j(c_i) f(\sigma_n(c_ih_n)) \,\doteq\, \gamma_j(\sigma_n), \qquad j=0,\dots,s-1,
\end{equation}
where $\,\doteq\,$ means {\em equal within machine precision}. Because of this fact, and for sake of brevity, we shall continue using $\sigma_n$ to denote the fully discrete approximation: 
\begin{equation}\label{sigman_d}
\sigma_n(ch_n) ~=~ \phi_{n-1}^{\aa,s}(c)~ +~ h_n^\aa\sum_{j=0}^{s-1} I^\aa P_j(c)\gamma_j^n, \qquad c\in[0,1].
\end{equation}
We observe that, in order to compute the (discrete) Fourier coefficients, it is enough to evaluate the memory term only at the quadrature abscissae $c_1,\dots,c_k$. In so doing, from (\ref{gammajn})-(\ref{sigman_d}) one obtains a discrete problem in the form:
\begin{equation}\label{dispro}
\gamma_j^n = \sum_{i=1}^k b_i P_j(c_i) f\left(\phi_{n-1}^{\aa,s}(c_i) + h_n^\aa\sum_{\ell=0}^{s-1} I^\aa P_\ell(c_i)\gamma_\ell^n\right), \qquad  j=0,\dots,s-1.
\end{equation}
Once it has been solved, by taking into account  (\ref{sigman_d}) and the fact that 
$I^\aa P_j(1)=\dfrac{\delta_{j0}}{\Gamma(\aa+1)}$, 
the approximation to $y(t_n)$ is given by:
\begin{equation}\label{ynew_d}
\bar y_n \,:=\, \sigma_n(h_n) \,=\,\phi_{n-1}^{\aa,s}(1) \,+\, \frac{h_n^\aa}{\Gamma(\aa+1)}\gamma_0^n.
\end{equation}
The equations  (\ref{dispro})--(\ref{ynew_d})) define a {\em Fractional HBVM with parameters $(k,s)$}. In short, FHBVM$(k,s)$. We refer to \cite{BBBI2023} for their analysis, and to \cite{BGI2024} for all implementation details.

\section{Numerical Tests}\label{num} 

In this section we report a few numerical tests using the Matlab\cpr code {\tt fhbvm}, implementing a FHBVM(22,20) method. The calling sequence of the code is:

\smallskip
\centerline{\tt [t,y,stats,err] = fhbvm( fun, y0, T, M )}

\no where:
\begin{itemize}
\item {\tt fun} is the function evaluating the r.h.s. of the equation (in vector mode), its Jacobian, and the order $\aa$ of the fractional derivative;
\item {\tt y0} is a matrix $\ell\times m$ containing the initial conditions in (\ref{ivp});
\item {\tt T} is the final integration time;
\item {\tt M} is an integer parameter such that $h_{\max} =~${\tt T/M} is approximately the largest stepsize in the final mesh ({\tt M} should be as small as possible);
\item {\tt t,y} contain the computed mesh and solution;
\item {\tt stats} (optional) is a vector containing some time statistics;
\item {\tt err} (optional) contains the estimate of the absolute error obtained on a doubled mesh. 
\end{itemize}
We refer to the inline help ({\tt help fhbvm}) for further details. The code is available at the url \cite{fhbvm}. All numerical tests have been done on a M2-Silicon based computer with 16GB of shared memory, using Matlab\cpr R2023b.
We shall also make a comparison with the Matlab\cpr code {\tt flmm2} \cite{Ga2018},\footnote{In particular, the BDF2 method is selected ({\tt method=3}), with the parameters {\tt tol=1e-15} and {\tt itmax=1000}.} in order to emphasize the potentialities of the new code. The comparisons will be done by using a so called {\em Work Precision Diagram (WPD)}, where the execution time (in {\tt sec}) is plotted against accuracy. The accuracy, in turn, is measured through the {\em mixed error significant computed digits} ({\tt mescd}) \cite{testset}, defined, by denoting $./$ the component-wise division, as\,\footnote{This definition corresponds to set {\tt atol=rtol} in the definition used in \cite{testset}.}
$${\tt mescd} := -\log_{10} \max_{i=0,\dots,N} \| (y(t_i)-\bar y_i)./(1+|y(t_i)|)\|_\infty,$$
with $t_i$, $i=0,\dots,N$, the computational mesh of the considered solver, $y(t_i)$ and $\bar y_i$ the corresponding values of the solution and of its approximation, and $|y(t_i)|$ the vector of the absolute values. The following parameters are used for all the examples below, to compute the corresponding WPD:
\begin{itemize}
\item {\tt flmm2}\,:\, $h=\frac{1}5 2^{-\nu}$, $\nu=1,\dots,20$; 

\smallskip
\item {\tt fhbvm}\,:\, $M=2,3,4,5$.
\end{itemize}

\begin{figure}[t]
\centering
\includegraphics[width=9cm]{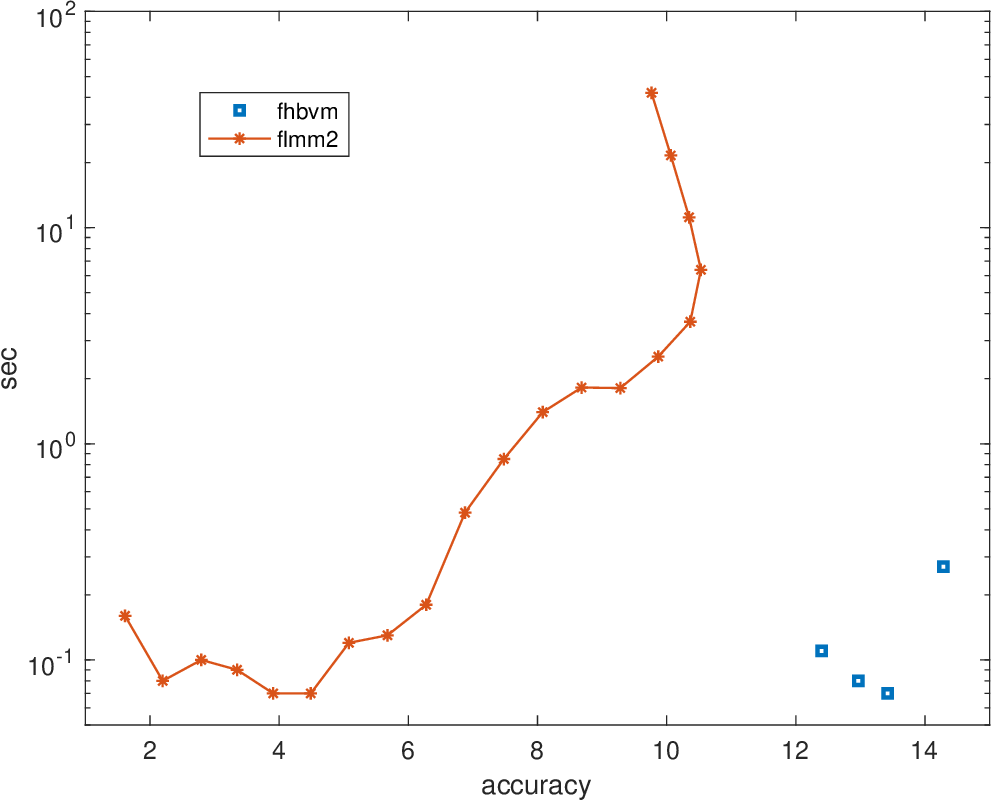}
\caption{Work-precision diagram for problem (\ref{ex1}), $\alpha=1.3$.}
\label{exe1}
\end{figure}

\paragraph{Example~1} 
The first problem is given by the Example~C.1, in \cite[page\,207]{Di2010}, for $\aa\in(1,2)$:
\begin{eqnarray}\nonumber
y^{(\aa)}(t) &=& -y(t)^{3/2} +\frac{40320}{\Gamma(9-\aa)}t^{8-\aa} -
           3\frac{\Gamma(5+\aa/2)}{\Gamma(5-\aa/2)}t^{4-\aa/2} +
          \left(\frac{3}2 t^{\aa/2} - t^4 \right)^3 + \frac{9}4\Gamma(\aa+1),\\
          && t\in[0,1], \qquad y(0)= y'(0)=0,\label{ex1}
\end{eqnarray} 
having unique solution 
$$ y(t) = t^8-3t^{4+\aa/2}+\frac{9}4t^\aa.$$
We consider the value $\aa=1.3$. Figure~\ref{exe1} contains the obtained WPD, from which one deduces that {\tt flmm2} can reach at most 10 digits of accuracy (further decreasing of the stepsize only increases the execution time, whereas the errors grow, due to round-off errors) in about 6 {\tt sec}. On the other hand, {\tt fhbvm} is able to gain full machine accuracy, with execution times of the order of one tenth of a second.

\begin{figure}[t]
\centering
\includegraphics[width=9cm]{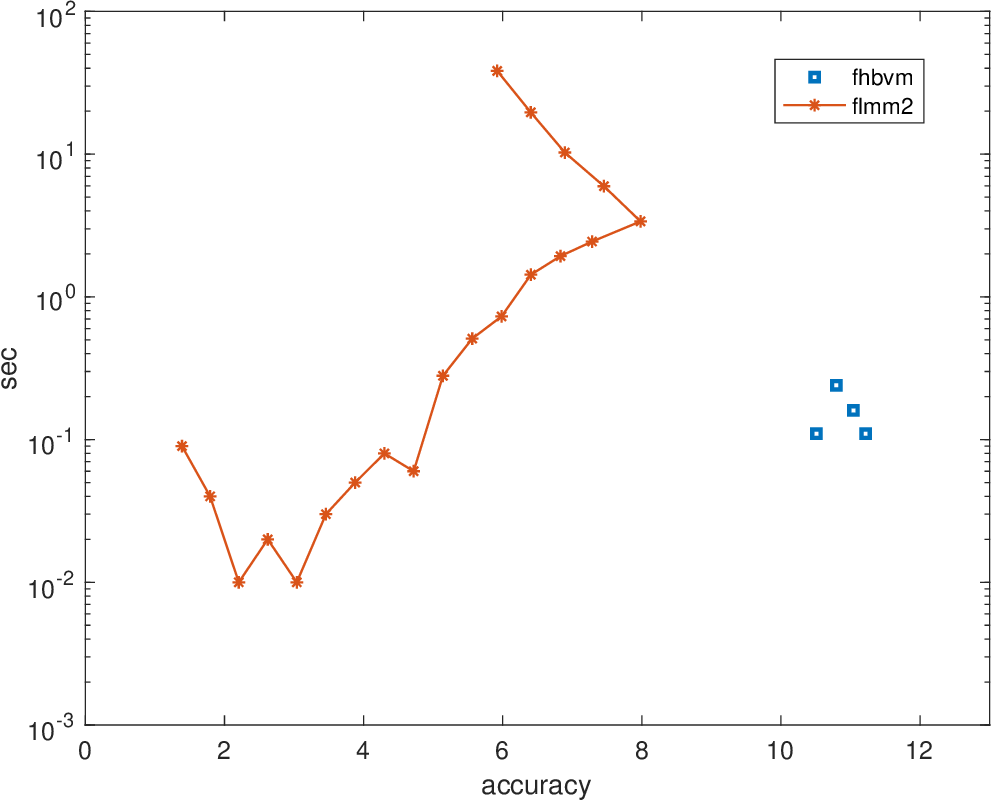}
\caption{Work-precision diagram for problem (\ref{ex2}).}
\label{exe2}
\end{figure}

\paragraph{Example~2} 
The next problem is taken from \cite[Example~2]{Sa2023}:
\begin{eqnarray}\nonumber
y^{(3/2)}(t)  &=& \frac{y(t)^2 -(t^{1.9}-1)^2}2 +\frac{\Gamma(2.9)}{\Gamma(1.4)}t^{0.4}, \qquad t\in[0,1],\\
&& y(0) = -1,\quad y'(0)=0,\label{ex2}
\end{eqnarray}
whose solution is~ $y(t) = t^{1.9}-1$. Figure~\ref{exe2} contains the obtained results, from which one may see that {\tt flmm2} can reach at most 8 significant digits, whereas {\tt fhbvm} obtains 11 significant digits with a much smaller time.

\begin{figure}[t]
\centering
\includegraphics[width=9cm]{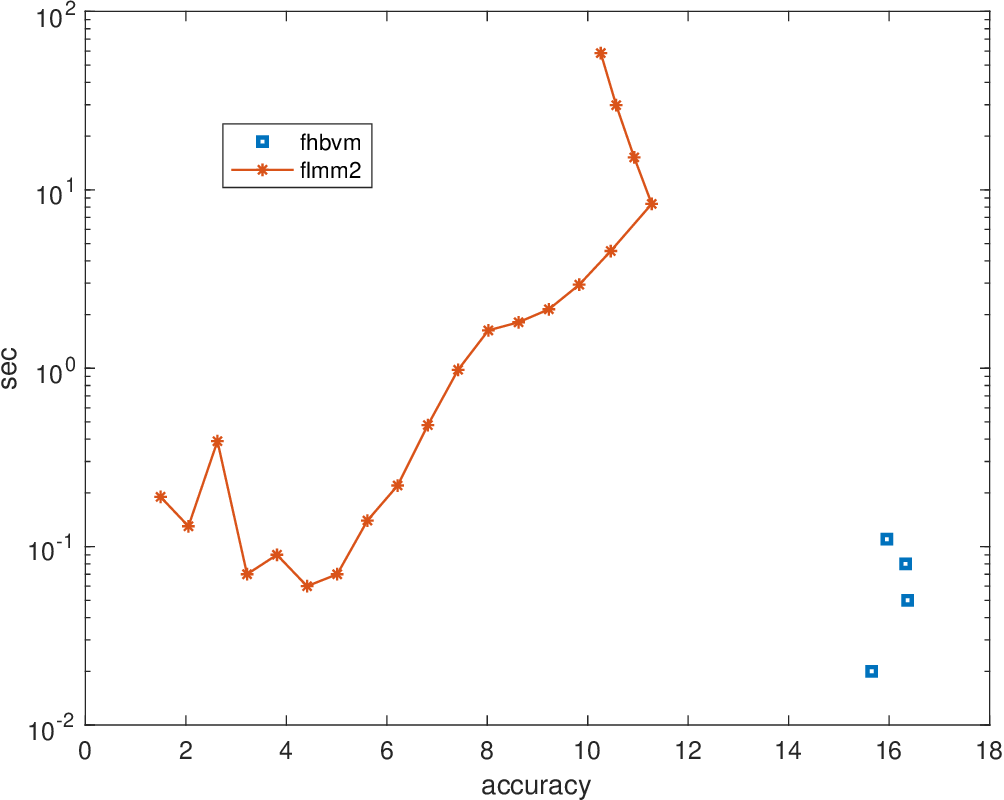}
\caption{Work-precision diagram for problem (\ref{ex3}), $\aa=1.25$.}
\label{exe3}
\end{figure}

\paragraph{Example~3} 
The third example is obtained by suitably combining Examples~2 and 3 in \cite{JN2017}:
 \begin{eqnarray}\nonumber
 y^{(\aa)}_1(t)&=&  \frac{\Gamma(4+\aa)}6t^3 -t^{8+2\aa}+y_2(t)^2,\\ \nonumber
 y^{(\aa)}_2(t)  &=& \frac{\Gamma(5+\aa)}{24}t^4 +t^{3+\aa}-y_1(t), \qquad t\in[0,1], \\
 && y_1(0) = y_1'(0) = y_2(0)=y_2'(0)=0, \label{ex3}
 \end{eqnarray}
 whose solution is given by
 $$y_1(t) = t^{3+\aa}, \qquad y_2(t) = t^{4+\aa}.$$
 We consider the same value $\aa=1.25$ considered in \cite{JN2017}. Figure~\ref{exe3} summarizes the obtained results, from which one may see that {\tt flmm2} reaches about 12 significant digits in about 8 {\tt sec}, while {\tt fhbvm} obtains about 17 significant digits in less than one tenth of a second. 
 
\begin{figure}[t]
\centering
\includegraphics[width=9cm]{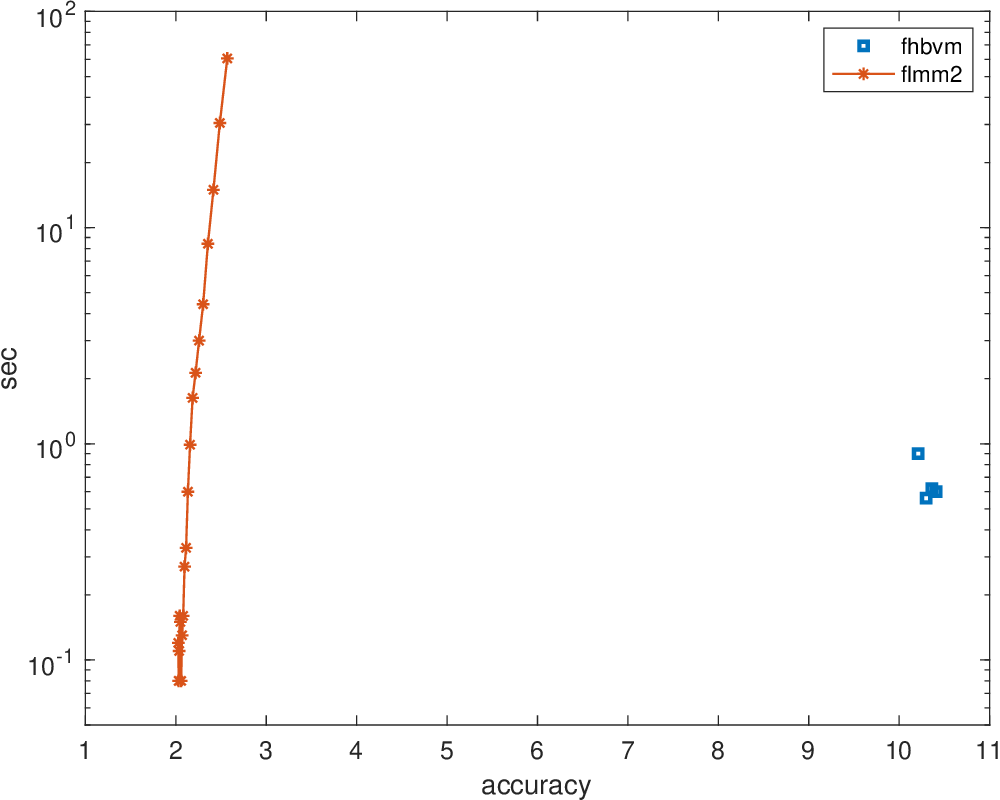} 
\caption{Work-precision diagram for problem (\ref{ex4}), $\aa=0.25$.}
\label{exe4}
\end{figure}

\paragraph{Example~4} 
At last, we consider the following {\em stiff} problem, which is adapted from Example~2 in \cite{BGI2024}:
 \begin{equation}\label{ex4}
 y^{(\aa)}(t)= \pmatrix{rr} -100 &0\\ -99 &-1\endpmatrix y(t), \qquad t\in[0,20], \qquad y(0) = \pmatrix{c}2\\ 3\endpmatrix. 
 \end{equation}
Its solution, for $\aa\in(0,1)$, is given by:
 $$y(t) = \pmatrix{c} 2 E_\aa(-100 t^\aa)\\ 2 E_\aa(-100 t^\aa)+E_\aa(-t^\aa)\endpmatrix,$$
 with $E_\aa$ the one-parameter Mittag-Leffler function. We consider the value $\aa=0.25$.
 Figure~\ref{exe4} summarizes the obtained results, from which one may see that {\tt flmm2}  reaches less than 3 significant digits in about  60 {\tt sec}, whereas {\tt fhbvm} reaches over 10 significant digits in less than 1 {\tt sec}. 
 
\section{Conclusions}\label{fine}

In this paper we have recalled the main facts about {\em Fractional HBVMs (FHBVMs)}, a class of numerical methods aimed at solving FDE-IVPs of Caputo type. In particular, the method FHBVM$(22,20)$ has been recently implemented in the Matlab\cpr code {\tt fhbvm} \cite{BGI2024}, and is freely available at the URL \cite{fhbvm}. Numerical comparisons with another state-of-art code confirm the potentialities of the approach and the efficiency of the code {\tt fhbvm}. In particular, the present version of the code has been extended to cope with high-order FDE-IVPs, w.r.t. the version presented in \cite{BGI2024}.
 
\paragraph*{\bf Acknowledgements.}\quad The authors are members of the Gruppo Nazionale Calcolo Scientifico-Istituto Nazionale di Alta Matematica (GNCS-INdAM).

\end{document}